\documentclass[a4paper,11pt,twoside,reqno]{amsart}

\usepackage[utf8]{inputenc}
\usepackage[plainpages=false,pdfpagelabels=true]{hyperref}
\usepackage{amssymb,amsthm}
\usepackage[margin=1in]{geometry}

\newtheorem{Satz}{Theorem}[section]
\newtheorem{Prop}[Satz]{Proposition}
\newtheorem{Lem}[Satz]{Lemma}
\newtheorem{Cor}[Satz]{Corollary}
\newcommand{\vol}{{\operatorname{vol}}}
\theoremstyle{definition}
\newtheorem{Dfn}[Satz]{Definition}
\newtheorem{Bem}[Satz]{Remark}
\newcommand{\zarg}{d\phi(e_1)\wedge d\phi(e_2)}
\newcommand{\uarg}{du(e_1)\wedge du(e_2)}
\newcommand{\uargt}{du_t(e_1)\wedge du_t(e_2)}
\newcommand{\Hom}{\operatorname{Hom}}
\newcommand{\hess}{\operatorname{Hess}}
\newcommand{\sff}{\mathrm{I\!I}}
\newcommand{\tr}{\operatorname{Tr}}

\allowdisplaybreaks[1]

\renewcommand{\epsilon}{\varepsilon}

\newcommand{\R}{\ensuremath{\mathbb{R}}}

\numberwithin{equation}{section}

\title{A global weak solution to the full bosonic string heat flow}
\author{Volker Branding}
\date{\today}
\address{University of Vienna, Faculty of Mathematics\\
Oskar-Morgenstern-Platz 1, 1090 Vienna, Austria\\}
\email{volker.branding@univie.ac.at}
\subjclass[2010]{58E20, 35K55, 53C80}
\keywords{full bosonic string; heat flow; global weak solution}
\begin{document}

\begin{abstract}
We prove the existence of a unique global weak solution to the full bosonic string
heat flow from closed Riemannian surfaces to an arbitrary target
under smallness conditions on the two-form and the scalar potential.
The solution is smooth with the exception of finitely many singular points.
Finally, we discuss the convergence of the heat flow and obtain a new existence result
for critical points of the full bosonic string action.
\end{abstract} 

\maketitle

\section{Introduction and Results}
The action functional for the full bosonic string is an important model
in contemporary theoretical physics. It is defined for a map
from a two-dimensional domain taking values in a manifold.
The action functional consists of three contributions:
Besides the \emph{Polyakov action} one considers the so-called \emph{B-field action} and a \emph{Dilaton}
contribution. For the physics background of the full bosonic string we refer to \cite[p. 108]{MR2151029}.

This article is a sequel to previous work concerning the existence of critical points of 
the full bosonic string action. 
In \cite{MR3573990} an existence result was given in the case of the domain being
a closed Riemannian surface and the target a Riemannian manifold
having negative sectional curvature. 

Moreover, a second existence result has been established in \cite{MR3624770}
for the domain being two-dimensional Minkowski space and the target an arbitrary
closed Riemannian manifold.

The aim of this article is to extend the existence result from \cite{MR3573990} to arbitrary
targets without posing any curvature assumption. 
In addition, we prove a regularity result for weak solutions of the critical points of
the full bosonic string action.

Let us explain the geometric setup in more detail.
Throughout this article \((M,h)\) is a closed Riemannian surface and \((N,g)\) a closed, oriented Riemannian manifold of dimension \(\dim N\geq 3\).
For a map \(\phi\colon M\to N\) we consider the square of its differential giving rise to the well-known Dirichlet energy,
whose critical points are \emph{harmonic maps}.
Let \(B\) be a two-form on \(N\), which we pull back by the map \(\phi\) and \(V\colon N\to\R\) be a scalar function.

In the physics literature the full action for the bosonic string is given by
\begin{align}
\label{energy-functional}
S_{bos}(\phi,h)=\int_M\big(\frac{1}{2}|d\phi|^2+\phi^\ast B+V(\phi)\big)d\vol_h.
\end{align}
We explicitly state the dependence of the action functional on the metric of the domain \(M\)
since the scalar potential \(V(\phi)\) is not invariant under conformal transformations.
Note that in the physics literature the scalar potential \(V(\phi)\) often gets multiplied with the scalar curvature of the domain.

In the mathematics literature there have been several articles dealing with energy functionals similar to \eqref{energy-functional}.
On the one hand there is the notion of \emph{harmonic maps with potential} introduced in \cite{MR1433176,MR1800592}, which are critical
points of \eqref{energy-functional} with \(B=0\).
On the other hand there have been several studies of the heatflow of \eqref{energy-functional} with \(V=0\), see for example
\cite{MR1393561,MR1424349} and \cite{MR1877837}.
For more references on the mathematical background see the introduction of \cite{MR3573990} and references therein.
The tools that we use in this article mostly originate from the theory of harmonic maps,
see \cite{MR826871} and the book \cite{MR1078018} for a detailed presentation.

The Euler-Lagrange equation of the functional \eqref{energy-functional} is given by
\begin{align}
\label{euler-lagrange-intrinsic}
\tau(\phi)=Z(\zarg)+\nabla V(\phi),
\end{align}
where \(\tau(\phi):=\tr\nabla d\phi\in\Gamma(\phi^\ast TN)\) denotes the tension field of the map \(\phi\) and the vector-bundle homomorphism \(Z \in \Gamma (\Hom(\Lambda^2T^\ast N,TN))\) is defined by the equation
\begin{align*}
\Omega(\eta,\xi_1,\xi_2)=\langle Z(\xi_1\wedge\xi_2),\eta\rangle,
\end{align*}
where \(\Omega=dB\) is a three-form on \(N\) and \(e_1,e_2\) an orthonormal basis of \(TM\).
For a derivation of \eqref{euler-lagrange-intrinsic} see \cite[Proposition 2.1]{MR3573990}.

First of all, we analyze the regularity of weak solutions of \eqref{euler-lagrange-intrinsic} and prove the following
\begin{Satz}
\label{theorem-regularity}
Let \((M,h)\) be a closed Riemannian surface and \((N,g)\) a closed Riemannian manifold with \(\dim N\geq 3\).
Suppose that \(\phi\in W^{1,2}(M,N)\) solves \eqref{euler-lagrange-intrinsic} in a distributional sense.
If \(V(\phi)\) is smooth then \(\phi\in C^\infty(M,N)\).
\end{Satz}

The major part of this article is devoted to the study of the \(L^2\)-gradient flow of the functional \eqref{energy-functional},
which is given by the following evolution equation
\begin{align}
\label{flow-intrinsic}
\frac{\partial\phi}{\partial t}(x,t)=&\tau(\phi)(x,t)-Z(\zarg)(x,t)-\nabla V(\phi)(x,t), \\
\nonumber\phi(x,0)=&\phi_0(x).
\end{align}
This is a natural generalization of the harmonic map heat flow from surfaces.
Although most of the analytical results obtained in this article follow the ideas from the standard
harmonic map heat flow we will encounter several new phenomena due to the presence of the scalar potential \(V(\phi)\)
in the action functional. For simplicity we will mostly assume the scalar potential is smooth. However,
we will point out the influence of a potential of lower regularity on the solution of \eqref{flow-intrinsic}
at several places.

By assumption the manifold \(N\) is compact, hence the potential \(V(\phi)\) satisfies
\(-A_1\leq V(\phi)\leq A_2\)
for positive constants \(A_1,A_2\). Exploiting this fact, we set
\begin{align}
\label{potential-shift}
0\leq\tilde V(\phi):=V(\phi)+A_1.
\end{align}

We will prove the following 
\begin{Satz}
\label{theorem-flow}
Let $(M,h)$ be a closed Riemannian surface and \((N,g)\) a closed Riemannian manifold.
Moreover, suppose that \(|B|_{L^\infty}<\frac{1}{2}\) and that \(V(\phi)\in C^\infty(N,\R)\) satisfies
\begin{align}
\label{assumption-smallness-potential}
\int_M\tilde V(\phi)d\vol_h<\delta
\end{align}
for some small \(\delta>0\).

Then for any initial data $\phi_0\in W^{1,2}(M,N)$ there exists a global weak solution 
\[
\phi\colon M\times [0,\infty)\to N
\]
of \eqref{flow-intrinsic} on \(M\times [0,\infty)\), which is smooth away from at most finitely many singular points $(x_k,t_k), 1\leq k\leq K$
with \(K=K(\phi_0,|V(\phi)|_{L^\infty},|B|_{L^\infty},M)\). 
The weak solution constructed here is unique and the energy functional \eqref{energy-functional} 
is decreasing with respect to time.\\
Moreover, there exists a sequence \(t_k\to\infty\) such that \(\phi(\cdot,t_k)\)
converges weakly in $W^{1,2}(M,N)$ to a solution of \eqref{euler-lagrange-intrinsic} denoted by
\(\phi_\infty\) as $k\to \infty$ suitably and strongly away from finitely many points \((x_k,t_k=\infty)\).
The limiting map \(\phi_\infty\) is smooth on \(M\setminus\{x_1,\ldots,x_K\}\).
\end{Satz}

Let us give a more precise definition of what is meant by a singularity in Theorem \ref{theorem-flow}.
We say that \((x_0,t_0)\) is a singular point of \eqref{flow-intrinsic} if for any \(R>0\)
\begin{align*}
\limsup_{t\to t_0}\int_{B_R(x_0)}|d\phi|^2d\mu\geq \delta_1,
\end{align*}
where \(\delta_1>0\) will be determined along the proof,
and \(B_R(x_0)\) denotes the geodesic ball around \(x_0\) with radius \(R\).

\begin{Bem}
\begin{enumerate}
\item Note that we can always perform a conformal rescaling of the metric \(h\) on the domain to achieve the smallness condition
 \eqref{assumption-smallness-potential}. Such a conformal transformation does not affect the other two terms in \eqref{energy-functional}
  since they are invariant under conformal transformations. 
\item In the case of the standard harmonic map heat flow from surfaces to general targets one
  can blow up the singular points that form along the flow. This procedures makes use of the fact
  that the harmonic map heat flow is invariant under parabolic rescaling. The inclusion of the scalar potential
  in the action functional \eqref{energy-functional} breaks the conformal invariance, as a consequence the critical
  points of \eqref{energy-functional} do not scale nicely. Hence, we cannot expect to blow up the singular points
  that form along \eqref{flow-intrinsic}.
 \item If we compare the results obtained in this article with the main results from \cite{MR3573990}
  we can make the following observations: In \cite{MR3573990} an existence result for \eqref{euler-lagrange-intrinsic}
  could be obtained under the assumption that the target manifold has negative curvature. In this article we do
  not impose any curvature condition on the target instead we have to make strong assumption on the scalar potential \(V(\phi)\).
\end{enumerate}
 
\end{Bem}

This article is organized as follows: In section 2 we study the regularity of weak solutions to \eqref{euler-lagrange-extrinsic}.
Afterwards, in section 3, we study the heat flow associated to \eqref{euler-lagrange-extrinsic} and prove Theorem \ref{theorem-flow}.

Whenever employing local coordinates, we will use Greek indices for coordinates on the domain and Latin indices 
for coordinates in the target. In addition, we will make use of the usual summation convention, that is we will
sum over repeated indices.

\section{Analytic aspects of the full bosonic string}
In this section we want to analyze several analytical properties of solutions of 
\eqref{euler-lagrange-intrinsic}. 
To this end we make use of the Nash embedding theorem 
and assume that \(N\subset\R^q\). Then \eqref{euler-lagrange-intrinsic} acquires the form
\begin{align}
\label{euler-lagrange-extrinsic}
\Delta u=\sff(du,du)+Z(du(e_1)\wedge du(e_2))+\nabla V(u),
\end{align}
where \(u\colon M\to\R^q\) and \(\sff\) denotes the second fundamental form of \(N\) in \(\R^q\).
For the equivalence of \eqref{euler-lagrange-intrinsic} and \eqref{euler-lagrange-extrinsic}
see \cite[Lemma 3.8]{MR3573990}.

In particular, we want to address the question how
the regularity of the scalar potential \(V(u)\) influences the regularity of
the solution of \eqref{euler-lagrange-intrinsic}. To this end, we will make the following
definition:
\begin{Dfn}
We call \(u\in W^{1,2}(M,N)\) a weak solution if it solves \eqref{euler-lagrange-extrinsic}
in a distributional sense.
\end{Dfn}

A similar study has already been performed in \cite{MR2415206,MR2250890} for 
harmonic maps with potential, that is critical points of \eqref{energy-functional} with \(B=0\).
Fortunately, by now there exist powerful tools that are well-adapted to \eqref{euler-lagrange-extrinsic}.
In the following we will make use of the following regularity result from \cite{MR3020100}.
\begin{Satz}
\label{topping-sharp}
Suppose that \(u\in W^{1,2}(D,\R^q)\) is a weak solution of
\begin{align}
-\Delta u=A\cdot\nabla u+f,\qquad f\in L^p(D,\R^q),
\end{align}
where \(A\in L^2(D,so(q)\otimes\R^2)\) and \(p\in(1,2)\). Then \(u\in W_{loc}^{2,p}(D)\).
In particular, if \(f=0\), then \(u\in W_{loc}^{2,p}\) for all \(p\in [1,2)\) and \(u\in W_{loc}^{1,q}\) for all \(q\in [1,\infty)\).
Moreover, for \(U\subset D\), there exist \(\eta_0=\eta_0(p,q)>0\) and \(C=C(p,m,U)<\infty\)
such that if \(\|A\|_{L^2(D)}\leq\eta_0\), then the following estimate holds
\begin{align}
\label{phi-w2p}
\|u\|_{W^{2,p}(U)}\leq C(\|f\|_{L^p(D)}+\|u\|_{L^1(D)}).
\end{align}
\end{Satz}

In order to be able to apply Theorem \ref{topping-sharp} we
need to rewrite the right hand side of \eqref{euler-lagrange-extrinsic}.
We denote coordinates in the ambient space \(\R^q\) by \((y^1,y^2,\ldots,y^{q})\).
Let \(\nu_l,l=n+1,\ldots,q\) be an orthonormal frame field for the normal bundle \(T^\perp N\).
For \(X,Y\in T_yN\) and \(\nabla_Y\nu_k=Y^i\frac{\partial\nu_k}{\partial y^i}\)
we express the second fundamental form as
\begin{equation*}
\sff_y(X,Y)=\langle X,\nabla_Y\nu_l\rangle\nu_l=X^iY^j\frac{\partial\nu^i_l}{\partial y^j}\nu_l.
\end{equation*}
Let \(D\) be a domain in \(M\) and consider a weak solution of \eqref{euler-lagrange-extrinsic}.
We choose local isothermal coordinates \(z=x+iy\), set \(e_1=\partial_x\), \(e_2=\partial_y\)
and use the notation \(u_\alpha=du(e_\alpha)\). 
Moreover, note that \(u_\alpha\in TN\) and \(\nu_l\in T^\perp N\), which implies that
\begin{align}
\label{orthogonality}
u^i_\alpha\nu^i_l=0
\end{align}
for all \(\alpha\). Hence, we may write
\begin{align}
\label{skew-sff}
\sff^m(u_\alpha,u_\alpha)=u^i_\alpha u^j_\alpha\big(\frac{\partial\nu_l^i}{\partial y^j}\nu_l^m-\frac{\partial\nu_l^m}{\partial y^j}\nu_l^i\big),\qquad m=1,\ldots,q, 
\end{align}
where we used \eqref{orthogonality} in the second term on the right hand side.
In addition, we note that
\[
Z^m(\uarg)=Z^m(\partial_{y^i}\wedge\partial_{y^j})u^i_xu^j_y,\qquad m=1,\ldots,q.
\]
By the definition of \(Z\) and exploiting the skew-symmetry of the three-form \(\Omega\), we find (see also \cite{MR3305429})
\begin{align}
\label{skew-z}
Z^k(\partial_{y^i}\wedge\partial_{y^j})=-Z^i(\partial_{y^k}\wedge\partial_{y^j}).
\end{align}
We are now in the position to show that solutions of \eqref{euler-lagrange-extrinsic}
have a structure such that Theorem \ref{topping-sharp} can be applied.

\begin{Prop}
Let \((M,h)\) be a closed Riemannian surface and let \((N,g)\) be a compact Riemannian manifold.
Assume that \(u\colon D\to N\) is a weak solution of \eqref{euler-lagrange-extrinsic}.
Let \(D\) be a simply connected domain of \(M\). 
Then there exists \(A^i_{~j}\in L^2(D,so(q)\otimes\R^2)\) such that
\begin{align}
\label{u-antisymmetric}
-\Delta u^m=A^m_{~i}\cdot\nabla u^i+(\nabla V(u))^m
\end{align}
holds.
\end{Prop}
\begin{proof}
By assumption \(N\subset\R^{q}\) is compact, we denote its unit normal field by \(\nu_l,l=n+1,\ldots,q\).
Using \eqref{skew-sff} and \eqref{skew-z}, we denote
\[
A^m_{~i}=\begin{pmatrix}
F^m_{~i} \\
G^m_{~i}
\end{pmatrix},
\qquad
i,m=1,\ldots,q
\]
with
\begin{align*}
F^m_{~i}:=&\left(\frac{\partial\nu_l^i}{\partial y^j}\nu_l^m-\frac{\partial\nu_l^m}{\partial y^j}\nu_l^i\right)u^j_x+Z^m(\partial_{y^i}\wedge\partial_{y^j})u^j_y,\\
G^m_{~i}:=&\left(\frac{\partial\nu_l^i}{\partial y^j}\nu_l^m-\frac{\partial\nu_l^m}{\partial y^j}\nu_l^i\right)u^j_y-Z^m(\partial_{y^i}\wedge\partial_{y^j})u^j_x.
\end{align*}
The skew-symmetry of \(A^m_{~i}\) can be read of from its definition and the properties of \(Z\), see \eqref{skew-z}.
By assumption \(u\) is a weak solution of \eqref{euler-lagrange-extrinsic}, hence \(A^m_{~i}\in L^2(D,so(q)\otimes\R^2)\)
completing the proof.
\end{proof}

First, we will assume that the scalar potential \(V(u)\) may have as little regularity
as possible. 

\begin{Cor}
Let \((M,h)\) be a closed Riemannian surface and let \(N\) be a compact Riemannian manifold.
Assume that \(u\colon D\to N\) is a weak solution of \eqref{euler-lagrange-extrinsic}.
Fix \(p\in (1,2)\) and assume that the scalar potential is of class \(V\in W^{1,p}(N,\R)\).
Then \(u\in W^{2,p}(M,N)\) and \(u\in W^{1,\frac{2p}{2-p}}(M,N)\).
\end{Cor}
\begin{proof}
This follows from Theorem \ref{topping-sharp} applied to \eqref{u-antisymmetric}
and the Sobolev embedding theorem in dimension two.
\end{proof}

One cannot expect to gain more regularity until one assumes that the potential \(V(u)\) 
has a better analytical structure.
In the case of a smooth potential \(V(u)\) we directly obtain Theorem \ref{theorem-regularity}.
\begin{proof}[Proof of Theorem \ref{theorem-regularity}]
This follows from elliptic regularity and a standard bootstrap argument.
\end{proof}

We conclude this section with the following ``gap-type'' theorem.

\begin{Prop}
Let \(u\) be a smooth solution of \eqref{euler-lagrange-extrinsic} with small energy
\(\|du\|_{L^2}<\epsilon\). Then the following inequality holds
\begin{align}
\|u\|_{W^{2,\frac{4}{3}}(M,N)}\leq C\|\nabla V\|_{L^\frac{4}{3}(M,N)},
\end{align}
where the positive constant \(C\) depends on \(M,N,\epsilon,|Z|_{L^\infty}\).
\end{Prop}
\begin{proof}
We estimate \eqref{euler-lagrange-extrinsic} as
\begin{align*}
\|\Delta u\|_{L^\frac{4}{3}(M,N)}\leq & C \||du|^2\|_{L^\frac{4}{3}(M,N)}+\|\nabla V\|_{L^\frac{4}{3}(M,N)} \\
\leq& C\|du\|_{L^2(M,N)}^2\|du\|_{L^4(M,N)}+\|\nabla V\|_{L^\frac{4}{3}(M,N)}.
\end{align*}
The claim follows by applying the Sobolev embedding theorem and choosing \(\epsilon\) sufficiently small.
\end{proof}
This allows us to draw the following
\begin{Cor}
If \(\|du\|_{L^2}\) is sufficiently small and 
\(\|\nabla V\|_{L^\frac{4}{3}(M,N)}\leq\delta \|u\|_{W^{2,\frac{4}{3}}(M,N)}\)
for \(\delta\) sufficiently small then \(u\) must be trivial.
\end{Cor}
Note that we do not have to make any assumption on \(V(u)\) but only on its gradient.

\section{The heat flow for the full bosonic string}
In this section we study the heat flow associated to \eqref{euler-lagrange-intrinsic}
and prove Theorem \ref{theorem-flow}. 

First, we will rewrite the action functional \eqref{energy-functional} in order to obtain
a functional that is easier to handle from an analytical point of view.

Shifting the potential \(V(\phi)\) as defined in \eqref{potential-shift}
we obtain the transformed energy functional
\begin{align}
\label{energy-functional-redefined-a}
\tilde S_{bos}(\phi,h)=\int_M\big(\frac{1}{2}|d\phi|^2+\phi^\ast B+\tilde V(\phi)\big)d\vol_h.
\end{align}
Note that \(\tilde S_{bos}(\phi,h)\geq 0\) if we also assume that \(|B|_{L^\infty}\leq\frac{1}{2}\).

\begin{Bem}
The critical points of \(\tilde S_{bos}(\phi,h)\) and \(S_{bos}(\phi,h)\)
coincide since both action functionals only differ by a constant. This fact is well known in physics:
The Lagrangian/Hamiltonian of a mechanical system can be changed by adding a constant since it does not contribute 
to the equations of motion. 
\end{Bem}

In order to deal with the analytic aspects of \eqref{flow-intrinsic} we again isometrically 
embedded the target manifold \(N\) into \(\R^q\). 
Then the corresponding heat-flow acquires the form
\begin{align}
\label{flow-rq}
\frac{\partial u_t}{\partial t}&=\Delta u_t-\sff(du_t,du_t)-Z(\uargt)-\nabla V(u_t), \\
\nonumber u(x,0)&=u_0(x),
\end{align} 
where \(u_t\colon M\times [0,T)\to\R^q\).
We will use a subscript \(t\) to denote the \(t\)-dependence of \(u\).
For a derivation of \eqref{flow-rq} see \cite[Lemma 4.1]{MR3573990}.
The existence of a shorttime solution can be obtained by standard methods, see for example \cite[Chapter 15]{MR2744149}.

\subsection{Energy Estimates}
In this subsection we will derive the necessary energy estimates for the study of \eqref{flow-rq}.

Let us introduce the following notation
\begin{align*}
E(u_t):=\int_M|du_t|^2d\vol_h, \\
E(u_t,B_{R}(x)):=\int_{B_{R(x)}}|du_t|^2d\mu.
\end{align*}
Here, \(B_R(x)\) denotes the geodesic ball of radius \(R\) around the point \(x\) and
by \(\iota_M\) we will denote the injectivity radius of \(M\).
Note that both of these energies are conformally invariant.

In addition, we introduce the following function space with \(Q=M\times [0,T)\) and \(dQ_h=d\vol_hdt\):
\begin{align*}
W:=\bigg\{\sup_{0\leq t\leq T}E(u_t)+\int_Q\big(|\nabla^2u_t|^2+\big|\frac{\partial u_t}{\partial t}\big|^2\big)dQ_h<\infty\bigg\}
\end{align*}

Due to the variational structure of our problem we have the following
\begin{Lem}
Let \(u_t\in W\) be a solution of \eqref{flow-rq}. Then the following equality holds
\begin{align}
\label{equality-energy}
\int_M(\frac{1}{2}|du_T|^2+u^\ast_T B+\tilde V(u_T))d\vol_h+\int_0^T\int_M\big|\frac{\partial u_t}{\partial t}\big|^2dQ_h
=\int_M(\frac{1}{2}|du_0|^2+u^\ast_0 B+\tilde V(u_0))d\vol_h.
\end{align}
\end{Lem}
\begin{proof}
We calculate
\begin{align*}
\frac{d}{dt}\frac{1}{2}\int_M|du_t|^2 d\vol_h=&-\int_M\langle\frac{\partial u_t}{\partial t},\Delta u_t\rangle d\vol_h \\
=&\int_M\big(-\big|\frac{\partial u_t}{\partial t}\big|^2-\langle\frac{\partial u_t}{\partial t},\sff(du_t,du_t)+Z(\uargt)+\nabla\tilde V(u_t)\rangle\big) d\vol_h \\
=&-\int_M\big|\frac{\partial u_t}{\partial t}\big|^2 d\vol_h-\int_M(\frac{\partial}{\partial t}u_t^\ast B+\frac{\partial}{\partial t}\tilde V(u_t))d\vol_h,
\end{align*}
where we used that \(\sff\perp\frac{\partial u_t}{\partial t}\).
The claim follows by integration with respect to \(t\).
\end{proof}

The next Lemma is the analogue of Lemma 3.6 from \cite{MR826871}.
\begin{Lem}
\label{lemma-energy-local}
Let \(u_t\in W\) be a solution of \eqref{flow-rq}.
For $R\in (0,i_M)$ and any $(x,t)\in Q$ there holds the estimate
\begin{align}
\label{local-energy-inequality}
\int_{B_R}(\frac{1}{2}|du_t|^2+u^\ast_t B+\tilde V(u_t))d\mu\leq\frac{C}{R^2}\int_Q|du_t|^2dQ_h+\int_{B_{2R}}(\frac{1}{2}|du_0|^2+u^\ast_0 B+\tilde V(u_0))d\mu,
\end{align}
where the constant $C$ only depends on $M$.
\end{Lem}
\begin{proof}
We choose a smooth cut-off function $\eta$ with the following properties
\begin{align*}
\eta\in C^\infty(M),\qquad\eta\geq 0,\qquad\eta=1~\textrm{on}~B_R(x_0),\\
\eta=0~\textrm{on}~M\setminus B_{2R}(x_0),\qquad |\nabla\eta|_{L^\infty}\leq\frac{C}{R},
\end{align*}
where again \(B_R(x_0)\) denotes the geodesic ball of radius \(R\) around \(x_0\in M\) and \(C\) a positive constant.
In addition, we choose an orthonormal basis \(\{e_\alpha,\alpha=1,2\}\) on \(M\)
such that \(\nabla_{e_\alpha}e_\beta=\nabla_{\partial_t}e_\alpha=0\) at the considered point.
By a direct calculation we obtain
\begin{align*}
\frac{\partial}{\partial t}\frac{1}{2}|du_t|^2=&d\langle\frac{\partial u_t}{\partial t},du_t\rangle
-\langle\frac{\partial u_t}{\partial t},\Delta u_t\rangle.
\end{align*}
Multiplying by the cut-off function $\eta^2$ and using the evolution equation \eqref{flow-rq}
we find
\begin{align*}
\frac{d}{dt}\frac{1}{2}\int_M\eta^2|du_t|^2d\vol_h
=&\int_M\big(\eta^2d\langle\frac{\partial u_t}{\partial t},du_t\rangle
+\eta^2\big(-\big|\frac{\partial u_t}{\partial t}\big|^2-\langle\frac{\partial u_t}{\partial t},Z(\uargt)\rangle \\
&-\langle\frac{\partial u_t}{\partial t},\nabla\tilde V(u_t)\rangle\big)d\vol_h \\
=&\int_M\big(\eta^2d\langle\frac{\partial u_t}{\partial t},du_t\rangle
-\eta^2\frac{\partial}{\partial t}u_t^\ast B-\eta^2\frac{\partial}{\partial t}\tilde V(u_t)-\eta^2\big|\frac{\partial u_t}{\partial t}\big|^2\big)d\vol_h.
\end{align*}
Using integration by parts we derive
\begin{align*}
\int_M\eta^2d\langle\frac{\partial u_t}{\partial t},du_t\rangle d\vol_h
\leq&2\int_M|\eta||d\eta||\frac{\partial u_t}{\partial t}||du_t|d\vol_h.
\end{align*}
Applying Young's inequality and by the properties of the cut-off function $\eta$, we find
\[
\frac{d}{dt}\int_M\eta^2(\frac{1}{2}|du_t|^2+u_t^\ast B+\tilde V(u_t))d\vol_h \leq\frac{C}{R^2}\int_M|du_t|^2d\vol_h.
\]
Integration with respect to $t$ yields the result.
\end{proof}

\begin{Prop}
Let \(u_t\in W\) be a solution of \eqref{flow-rq}. Moreover, 
suppose that \(|B|_{L^\infty}<\frac{1}{2}\).

Then the following monotonicity formulas hold
\begin{align}
\label{monotonicity-energy-global}
E(u_t)\leq &\delta_2\tilde S_{bos}(u_0,h) \\
\nonumber\leq&\delta_3 E(u_0)+\delta_2\int_M\tilde V(u_0)d\vol_h, \\
\label{monotonicity-energy-local}
E(u_t,B_R)\leq &C\delta_2^2\frac{T}{R^2}\tilde S_{bos}(u_0,h)
+\delta_2\tilde S_{bos}(u_0,B_{2R}) \\
\nonumber\leq &\delta_3 E(u_0,B_{2R})+C\delta_2^2\frac{T}{R^2}\tilde S_{bos}(u_0,h)
+\delta_2\int_{B_{2R}}\tilde V(u_0)d\mu,
\end{align}
where 
\begin{align*}
\delta_2:=\frac{1}{\frac{1}{2}-|B|_{L^\infty}},\quad\delta_3:=\frac{\frac{1}{2}+|B|_{L^\infty}}{\frac{1}{2}-|B|_{L^\infty}}.
\end{align*}
Here, \(S_{bos}(u_0,B_{2R})\) denotes the action functional at time \(0\) restricted to the ball \(B_{2R}\).
\end{Prop}

\begin{proof}
This follows from combining \eqref{equality-energy} and \eqref{local-energy-inequality}
and making use of the assumptions.
\end{proof}

In the following we want to control the energy of \(u_t\) locally.
\begin{Lem}
Let \(u_t\in W\) be a solution of \eqref{flow-rq}. Moreover, 
suppose that \(|B|_{L^\infty}<\frac{1}{2}\).
Then for any \(\delta_1>0\) there exist $R_1\in(0,i_M)$ and \(T_1>0\) such that
\begin{align}
\label{small-delta-one}
\sup_{\genfrac{}{}{0pt}{}{x\in M}{0\leq t\leq T_1}} E(u_t,B_{R_1})<\delta_1.
\end{align}
\end{Lem}
\begin{proof}
Given any \(u_0\) we can always find some \(R_1>0\) such that
\begin{align*}
\tilde S_{bos}(u_0,B_{2R_1})<\frac{\delta_1}{2\delta_2}
\end{align*}
for a positive constant $\delta_1$.
The statement then follows from \eqref{monotonicity-energy-local} by choosing
\begin{align*}
T_1=\frac{\delta_1}{2}\frac{R_1^2}{C\delta_2^2\tilde S_{bos}(u_0,h)}.
\end{align*}
\end{proof}

Let \(X\in\R^2\) be a bounded domain. Then Ladyzhenskaya's inequality holds, that is
\begin{Lem}
Assume that \(v\in W^{1,2}(X)\). Then the following inequality holds:
\begin{align*}
\|v\|^4_{L^4(X)}\leq C\|v\|^2_{L^2(X)}\|dv\|^2_{L^2(X)}
\end{align*}
\end{Lem}

In the following we need a local version of Ladyzhenskaya's inequality from above.
\begin{Lem}
Assume that \(v\in W\). Then there exists a constant $C$ such that for any \(R\in(0,i_M)\)
the following inequality holds:
\begin{align}
\label{local-sobolev-inequality}
\int_M|dv|^4 d\vol_h\leq C\sup_{x\in M}\int_{B_{R}(x)\author{}}|dv|^2d\vol_h\big(\int_M|\nabla^2v|^2d\vol_h+\frac{1}{R^2}\int_M|dv|^2d\vol_h\big).
\end{align}
\end{Lem}
\begin{proof}
A proof can for example be found in \cite[Lemma 6.7]{MR2431434}. 
\end{proof}

Making use of the Ricci identity we obtain the following formula for \(v\colon M\to\R^q\)
\begin{align}
\label{ricci-identity-second-derivatives}
\int_M|\Delta v|^2d\vol_h=\int_M|\nabla^2v|^2d\vol_h -\frac{1}{2}\int_M\operatorname{Scal}|dv|^2d\vol_h.
\end{align}
Here, \(\operatorname{Scal}\) denotes the scalar curvature of \(M\).

As a next step we will control the \(L^2\)-norm of the second derivatives of \(u_t\).
\begin{Lem}
Let \(u_t\in W\) be a solution of \eqref{flow-rq} and
suppose that \eqref{small-delta-one} holds with \(\delta_1\) sufficiently small.
Moreover, suppose that \(|B|_{L^\infty}<\frac{1}{2}\).
Then the following inequality holds
\begin{align}
\int_Q|\nabla^2u_t|^2dQ_h\leq C(1+\frac{T}{R^2}),
\end{align}
where the constant \(C\) depends on \(M,N,\delta_1,|Z|_{L^\infty},|B|_{L^\infty},|\hess V|_{L^\infty}\).
\end{Lem}
\begin{proof}
By a direct calculation we find
\begin{align*}
\frac{d}{dt}\frac{1}{2}\int_M\big|du_t|^2&d\vol_h=-\int_M\langle\Delta u_t,\frac{\partial u_t}{\partial t}\rangle d\vol_h \\
=&\int_M(-|\Delta u_t|^2+\langle\Delta u_t,\sff(du_t,du_t)+Z(\uargt)+\langle\Delta u_t,\nabla V(u_t)\rangle)d\vol_h\\
\leq &-\frac{1}{2}\int_M|\Delta u_t|^2d\vol_h +C\int_M|du_t|^4d\vol_h-\int_M\hess(du_t,du_t)d\vol_h.
\end{align*}
By assumption \(N\) is compact and we can estimate the Hessian of \(V\) by its maximum.  
Making use of \eqref{local-sobolev-inequality} and \eqref{ricci-identity-second-derivatives} we obtain
\begin{align*}
\frac{d}{dt}\frac{1}{2}\int_M\big|du_t|^2d\vol_h=&-\frac{1}{2}\int_M|\nabla^2u_t|^2d\vol_h+C\int_M|du_t|^2d\vol_h \\
&+C\sup_{x\in M}\int_{B_{R}(x)\author{}}|du_t|^2d\vol_h\big(\int_M|\nabla^2u_t|^2d\vol_h+\frac{1}{R^2}\int_M|du_t|^2d\vol_h\big).
\end{align*}
Choosing \(\delta_1\) small enough we get the following inequality
\begin{align*}
\frac{d}{dt}\frac{1}{2}\int_M\big|du_t|^2d\vol_h+C\int_M|\nabla^2u_t|^2d\vol_h\leq 
C\int_M|du_t|^2d\vol_h+\frac{C}{R^2}\int_M|du_t|^2d\vol_h.
\end{align*}
The claim follows by integration with respect to \(t\).
\end{proof}

Using the bound on the second derivatives, we can apply the Sobolev embedding theorem
to bound $\int_Q|du_t|^4dQ_h$.
\begin{Cor}
Let \(u_t\in W\) be a solution of \eqref{flow-rq} with \(\delta_1\) sufficiently small.
Moreover, suppose that \(|B|_{L^\infty}<\frac{1}{2}\).
Then we have for all $t\in[0,T_1)$
\begin{align}
\label{l4-small}
\int_Q|du_t|^4dQ_h\leq Cf(T_1)
\end{align}
with \(f(T_1)\) satisfying \(f(T_1)\to 0\) as \(T_1\to 0\).
\end{Cor}
\begin{proof}
The bound follow from \eqref{local-sobolev-inequality} and the previous estimate.
\end{proof}

As a next step we control the \(L^2\)-norm of the derivatives of \(u_t\) with respect to \(t\).
\begin{Lem}
Let \(u_t\in W\) be a solution of \eqref{flow-rq}. 
Moreover, suppose that \(|B|_{L^\infty}<\frac{1}{2}\).
If 
\(\sup_{(x,t)\in M\times[0,T_1)}E(u_t,B_{R_1}(x)))<\delta_1\) is small enough, we find for \(\xi>0\)
\begin{align}
\label{l2-partial-u-t}
\sup_{2\xi\leq t\leq T_1}\int_M\big|\frac{\partial u}{\partial t}(\cdot,t)|^2d\vol_h\leq C(1+\xi^{-1}),
\end{align}
where the positive constant \(C\) depends on \(M,N,\delta_1,u_0,|B|_{L^\infty},|Z|_{L^\infty},|\nabla Z|_{L^\infty},|\hess V|_{L^\infty}\).
\end{Lem}
\begin{proof}
By a direct calculation using \eqref{flow-rq} we find
\begin{align*}
\frac{d}{dt}\frac{1}{2}\int_M\big|\frac{\partial u_t}{\partial t}\big|^2d\vol_h=&-\int_M\big|\nabla\frac{\partial u_t}{\partial t}\big|^2d\vol_h
-\int_M\langle\frac{\nabla}{\partial t}\big(\sff(du_t,du_t)\big),\frac{\partial u_t}{\partial t}\rangle d\vol_h \\
&-\int_M\big(\langle\frac{\nabla}{\partial t}\big(Z(\uargt)\big),\frac{\partial u_t}{\partial t}\rangle 
-\hess V(\frac{\partial u_t}{\partial t},\frac{\partial u_t}{\partial t})\big)d\vol_h.
\end{align*}
Again, we can estimate the Hessian of the potential \(V(u_t)\) by its maximum since \(N\) is compact.
Consequently, we obtain
\begin{align*}
\frac{d}{dt}\frac{1}{2}\int_M\big|\frac{\partial u_t}{\partial t}\big|^2d\vol_h\leq&-\int_M\big|\nabla\frac{\partial u_t}{\partial t}\big|^2d\vol_h \\
&+C\int_M(|du_t|^2\big|\frac{\partial u_t}{\partial t}\big|^2+|du_t|\big|\nabla\frac{\partial u_t}{\partial t}\big|\big|\frac{\partial u_t}{\partial t}\big|
+\big|\frac{\partial u_t}{\partial t}\big|^2)d\vol_h \\
\leq &-\frac{1}{2}\int_M\big|\nabla\frac{\partial u_t}{\partial t}\big|^2d\vol_h+C\int_M|du_t|^2\big|\frac{\partial u_t}{\partial t}\big|^2 d\vol_h
+C\int_M\big|\frac{\partial u_t}{\partial t}\big|^2 d\vol_h.
\end{align*}
To control the second term on the right hand side we use another type of Sobolev inequality
(similar to \eqref{local-sobolev-inequality} for \(|t-s|\leq 1\)),
that is
\begin{align*}
\int_s^t\int_M|du_t|^2|\frac{\partial u_t}{\partial t}|^2&dQ_h\\ \leq
&\big(\int_s^t\int_M|du_t|^4dQ_h\big)^\frac{1}{2}
\big(\sup_{s\leq\theta\leq t} \int_M|\frac{\partial u}{\partial t}(\cdot,\theta)|^2 d\vol_h
+\int_s^t\int_M |\nabla\frac{\partial u_t}{\partial t}|^2dQ_h\big).
\end{align*}
Using \eqref{l4-small} and 
integrating over a small time interval $t-s<z$,
we can absorb part of the right hand side in the left and obtain
\begin{align*}
\int_M|\frac{\partial u}{\partial t}(\cdot,t)|^2d\vol_h
\leq \inf_{t-z\leq s\leq t}C\int_M|\frac{\partial u}{\partial t}(\cdot,s)|^2d\vol_h+\delta_2\tilde S_{bos}(u_0,h).
\end{align*}
Finally, we estimate the infimum by the mean value, more precisely
\begin{align*}
\sup_{2\xi\leq t\leq T_1}\int_M|\frac{\partial u}{\partial t}(\cdot,t)|^2d\vol_h
\leq C(1+\xi^{-1})\int_s^t\int_M|\frac{\partial u_t}{\partial t}|^2dQ_h+C
\leq C(1+\xi^{-1}).
\end{align*}
Hence, we get the desired bound.
\end{proof}

\begin{Lem}
Let \(u_t\in W\) be a solution of \eqref{flow-rq} with \(|B|_{L^\infty}<\frac{1}{2}\).
As long as \(\delta_1\) is sufficiently small 
we have the following bound
\begin{align}
\label{bound-second-derivatives}
\int_M|\nabla^2u_t|^2d\vol_h\leq C.
\end{align}
The constant \(C\) depends on  \(M,N,\delta_1,u_0,|B|_{L^\infty},|Z|_{L^\infty},|\nabla Z|_{L^\infty},|\nabla V|_{L^\infty},|\hess V|_{L^\infty}\).
\end{Lem}
\begin{proof}
Using \eqref{flow-rq} and \eqref{ricci-identity-second-derivatives} we obtain the following inequality 
\begin{align*}
\int_M|\nabla^2u_t|^2d\vol_h\leq & \int_M|\Delta u_t|^2d\vol_h+C\int_M|du_t|^2d\vol_h \\
\leq & C\int_M\big(\big|\frac{\partial u_t}{\partial t}\big|^2+|du_t|^4+|\nabla V(u_t)|^2+|du_t|^2\big)d\vol_h.
\end{align*}
Since we are assuming the potential \(V(u)\) to be smooth and \(N\) to be compact we can easily estimate \(|\nabla V(u)|^2\).
Applying \eqref{local-sobolev-inequality} and \eqref{l2-partial-u-t}
with \(\delta_1\) sufficiently small yields the claim.
\end{proof}

\begin{Prop}[Higher Regularity]
Let \(u_t\in W\) be a solution of \eqref{flow-rq}. 
As long as \(\delta_1\) is small enough the solution \(u_t\) of \eqref{flow-rq}
is smooth.
\end{Prop}
\begin{proof}
This follows from standard regularity theory arguments, see for example \cite[Lemma 6.11]{MR2431434}
and references therein for more details.
\end{proof}

Let us close this section with the following remarks.
\begin{Bem}
\begin{enumerate}
 \item The fact that \(u_t\) is smooth as long as \(\delta_1\) is sufficiently small relies on the fact that
  we have a smooth scalar potential \(V(u)\). If we would assume lower regularity of \(V(u)\) then \(u_t\) would
  also have less regularity. We can use the parabolicity of \eqref{flow-rq} to smoothen out distributional initial data,
  but the parabolicity cannot compensate for a potential of bad regularity. 
  
  In order to achieve that \(u\in W^{2,2}(M,N)\) we have to require that \(V\in C^2(N,\R)\).
  By the Sobolev embedding theorem we then get that \(u\) is continuous, to gain more regularity we need better
  regularity of \(V(u)\).
 \item In the case of a smooth heat flow one can also consider the case of a non-compact target \(N\)
  and use the potential \(V(u)\) to constrain the image of \(M\) under \(u_t\) to a compact set.
  However, this argument makes use of the maximum principle, which we cannot apply in our case.
\end{enumerate}
\end{Bem}

\subsection{Longtime Existence}
In this section we establish the existence of a unique global weak solution to \eqref{flow-rq}
for all times \(t\in[0,\infty)\). Moreover, we will show that only finitely
many singularities will occur along the flow. First, we will give a uniqueness result.
\begin{Prop}
\label{prop-uniqueness}
Let \(u_t,v_t\in W\) be two solutions of \eqref{flow-rq} and suppose that \(|B|_{L^\infty}<\frac{1}{2}\).
If their initial data coincides, that is
\(u_0=v_0\), then \(u_t=v_t\) for all \(t\in[0,T)\).
\end{Prop}

\begin{proof}
Throughout the proof \(C\) will denote a universal constant that may change from line to line.
Let \(u_t,v_t\) be two solutions of \eqref{flow-rq}. We set \(w_t:=u_t-v_t\). 
By projecting to a tubular neighborhood \(\sff(u_t)(du_t,du_t), Z(u_t)(du_t(e_1)\wedge du_t(e_2))\) and \(\nabla V(u_t)\) can be thought of as vector-valued functions in \(\R^q\),
for more details see \cite[Lemma 4.8]{MR3573990}.
Exploiting this fact a direct computation yields
\begin{align*}
\frac{\partial w_t}{\partial t}=&\Delta w_t+\langle\sff(u_t)(du_t,du_t)-\sff(v_t)(dv_t,dv_t),w_t\rangle \\
&+\langle Z(u_t)(du_t(e_1)\wedge du_t(e_2))-Z(v_t)(dv_t(e_1)\wedge dv_t(e_2)),w_t\rangle
+\langle\nabla V(u_t)-\nabla V(v_t),w_t\rangle.
\end{align*}
Rewriting
\begin{align*}
\sff(u_t)(du_t,du_t)-\sff(v_t)(dv_t,dv_t)=&(\sff(u_t)-\sff(v_t))(du_t,du_t) \\
&+\sff(v_t)(du_t-dv_t,du_t)+\sff(v_t)(dv_t,du_t-dv_t) 
\end{align*}
and similarly for the terms containing \(Z\) we find
\begin{align*}
\frac{d}{dt}\frac{1}{2}\int_M|w_t|^2d\vol_h\leq&-\int_M|dw_t|^2d\vol_h
+C\int_M(|w_t|^2|du_t|^2+|w_t|^2|dv_t|^2)d\vol_h \\
&+C\int_M(|w_t||dw_t||du_t|+|w_t||dw_t||dv_t|)d\vol_h \\
&+\int_M|\langle\nabla V(u_t)-\nabla V(v_t),w_t\rangle| d\vol_h.
\end{align*}
This leads to the following inequality
\begin{align*}
\frac{1}{2}\|w_t\|^2_{L^2(M)}+\|dw_t\|^2_{L^2(Q)}
\leq &C\big(\|w_t\|^2_{L^4(Q)}\|dv_t\|^2_{L^4(Q)}
+\|w_t\|^2_{L^4(Q)}\|du_t\|^2_{L^4(Q)} \\
&+\|w_t\|_{L^4(Q)}\|du_t\|_{L^4(Q)}\|dw_t\|_{L^2(Q)} \\
&+\|w_t\|_{L^4(Q)}\|dv_t\|_{L^4(Q)}\|dw_t\|_{L^2(Q)} \\
&+\||\nabla V(u_t)-\nabla V(v_t)||w_t|\|_{L^1(Q)}\big).
\end{align*}
By assumption the scalar potential \(V(u)\) is smooth such that
we can apply the mean-value theorem and estimate
\begin{align*}
\||\nabla V(u_t)-\nabla V(v_t)||w_t|\|_{L^1(Q)}\leq C\|w_t\|^2_{L^2(Q)}.
\end{align*}
Using \eqref{l4-small} we obtain
\begin{align*}
\frac{1}{2}\int_M|w_t|^2d\vol_h+&\frac{1}{2}\int_0^T\int_M|dw_t|^2dQ_h \\
\leq& Cf(T)\big(\int_Q|w_t|^4dQ_h\big)^\frac{1}{2}+C\int_Q|w_t|^2dQ_h \\
\leq& Cf(T)\big(\sup_{[0,T]}\int_M|w_t|^2d\vol_h+\int_0^T\int_M|dw_t|^2dQ_h\big)+C\int_Q|w_t|^2dQ_h
\end{align*}
with \(f(T)\to 0\) as \(T\to 0\).
Taking the limit \(T\to 0\) and applying the Gronwall inequality allows us to conclude the claim.
\end{proof}

\begin{Bem}
The proof of the previous Proposition requires the potential \(V(\phi)\) to be sufficiently regular
such that we can apply the mean-value theorem. If we would allow for a potential with less regularity
it does not seem to be possible to prove uniqueness of the solution of \eqref{flow-rq}.
\end{Bem}

By the same strategy as in the case of standard harmonic maps we can establish 
the long-time existence of \eqref{flow-rq}.
\begin{Prop}[Long-time Existence]
Let \(u_t\in W\) be a solution of \eqref{flow-rq}.
Moreover, suppose that \(|B|_{L^\infty}<\frac{1}{2}\).
Then \eqref{flow-rq} admits a unique weak solution for $0\leq t<\infty$.
\end{Prop}
\begin{proof}
The first singular time \(T_0\) is characterized by the condition
\[
\limsup_{t\to T_0}E(u_t,B_R(x))\geq \delta_1.
\]
Since we have $\partial_tu_t \in L^2(M\times [0,T_0))$ and
also $E(u_t)\leq \delta_2 \tilde S_{bos}(u_0,h)$ for $0<t<T_0$, there exists 
\[u(\cdot,T_0)\in W^{1,2}(M,N)
\]
such that 
\[u(\cdot,t)\to u(\cdot,T_0)\]
weakly in \(W^{1,2}(M,N)\) as $t$ approaches $T_0$. 
In particular, we have
\[
E(u_{T_0})\leq\liminf_{s\to T_0} E(u_s)\leq \delta_2 \tilde S_{bos}(u_0,h),\qquad 0\leq t\leq T_0.
\]
Let
\(\tilde{u}_t\colon M\times[T_0,T_0+T_1)\to N\)
be a solution of \eqref{flow-rq}.
Assume that $\tilde{u}(x,t)=u(x,t)$.
We define
\[
\hat{u}_t=
\begin{cases}
u, &\qquad 0\leq t\leq T_0,\\
\tilde{u}_t, &\qquad T_0\leq t\leq T_0+T_1.
\end{cases}
\]
Now
\(
\hat{u}_t\colon M\times[0,T_0+T_1)\to N
\)
is a weak solution of \eqref{flow-rq}.
By iteration, we obtain a weak solution $u_t$ on a maximal time interval $T_0+\delta$
for some $\delta>0$. If $T_0+\delta<\infty$ the above argument shows that the solution $u_t$
may be extended to infinity, hence $T_0+\delta=\infty$. 
The uniqueness of the solution follows from Proposition \ref{prop-uniqueness}.
\end{proof}

\begin{Prop}
Let \(u_t\in W\) be a solution of \eqref{flow-rq}.
Suppose that \(|B|_{L^\infty}<\frac{1}{2}\) and that the scalar potential \(V(u)\) is
sufficiently small, that is
\begin{align}
\label{smallness-scalar-potential}
\int_M\tilde V(u_t)d\vol_h\leq \frac{\delta_1}{\delta_2}.
\end{align}
Then there are only finitely many singular points \((x_k,t_k),1\leq k\leq K\).
The number \(K\) depends on \(M,|V(u)|_{L^\infty},|B|_{L^\infty},u_0\).
\end{Prop}
\begin{proof}
We follow the presentation in \cite[p.138]{MR2431658} for the harmonic map heat flow.
We assume that $T_0>0$ is the first singular time and define the singular set as
\begin{align*}
Z(u,T_0)=\bigcap_{R>0}\big\{
x\in M\mid\limsup_{t\to T_0}E(u_t,B_R(x))\geq \delta_1
\big\}.
\end{align*}
Now, let $\{x_j\}^K_{j=1}$ be any finite subset of $Z(u,T_0)$.
Then we have for $R>0$
\[
\limsup_{t\to T_0}\int_{B_R(x_j)}|du_t|^2d\mu\geq\delta_1,\qquad 1\leq j\leq K.
\]
We choose $R>0$ such that all the $B_{2R}(x_j),1\leq j\leq K$ are mutually disjoint.
Then, we have by \eqref{monotonicity-energy-local}
\begin{align*}
K\delta_1\leq& \sum_{j=1}^K\limsup_{t\to T_0} E(u_t,B_R(x_j)) \\
\leq& \sum_{j=1}^K\big(\delta_2\limsup_{t\to T_0} \tilde S_{bos}(u_\xi,B_{2R}(x_j))+\frac{\delta_1}{2}\big) \\
\leq& \delta_2 \tilde S_{bos}(u_\xi,h)+\frac{K\delta_1}{2} \\
\leq& \delta_2 \tilde S_{bos}(u_0,h)+\frac{K\delta_1}{2}
\end{align*}
for any $\xi\in[T_0-\frac{R^2}{C\delta_2^2 \tilde S_{bos}(u_0,h)},T_0]$.
We conclude that
\[
K\leq \frac{2\delta_2}{\delta_1} \tilde S_{bos}(u_0,h),
\]
which implies the finiteness of the singular set \(Z(u,T_0)\).
Next, we show that there are only finitely many singular spatial points. We set
\[
\tilde{M}=M\setminus\bigcup_{1\leq j\leq K}B_{2R}(x_j)
\]
and calculate
\begin{align}
\label{inequality-F-singularities}
E(u_{T_0})=&\lim_{R\to 0} E(u_{T_0},\tilde{M})\\
\nonumber\leq&\lim_{R\to 0}\limsup_{t\to T_0} E(u_t,\tilde{M})\\
\nonumber =&\limsup_{t\to T_0}E(u_t)-\lim_{R\to 0}\sum_{j=1}^K\liminf_{t\to T_0} E(u_t,B_{2R}(x_j))\\
\nonumber\leq&\delta_2 \tilde S_{bos}(u_0,h)-\lim_{R\to 0}\sum_{j=1}^K\limsup_{t\to T_0} E(u_t,B_{R}(x_j))\\
\nonumber\leq&\delta_2 \tilde S_{bos}(u_0,h)-K\delta_1 \\
\nonumber\leq&\delta_3 E(u_0)+\delta_2\int_M\tilde V(u_0)d\vol_h-K\delta_1.
\end{align}
Suppose $T_0<\ldots<T_j$ are $j$ singular times and by $K_0,\ldots,K_j$
we denote the number of singular points at each singular time.
Set
\[
u_i=\lim_{t\to T_i}u_t,\qquad V_i=\lim_{t\to T_i}\int_M\tilde V(u_t)d\vol_h \qquad 0\leq i\leq j.
\]
By iterating \eqref{inequality-F-singularities} we get
\begin{align*}
E(u_j)\leq&\delta_3 E(u_{j-1})+\delta_2V_{j-1}-\delta_1 K_{j-1} \\
\leq& \delta_3^2 E(u_{j-2})-\delta_1(K_{j-1}+\delta_3 K_{j-2})+\delta_2(V_{j-1}+\delta_3V_{j-2})\\
\leq&\ldots \\
\leq&\delta_3^jE(u_0)+\sum_{i=0}^{j-1}\delta_3^{j-i-1}(\delta_2 V_i-\delta_1 K_i),
\end{align*}
which can be rearranged as
\begin{align*}
\sum_{i=0}^{j-1}\delta_3^{-i-1}(\delta_1K_i-\delta_2V_i)\leq E(u_0).
\end{align*}
In order to conclude that the number of singularities is finite we have to ensure that
\begin{align*}
\delta_1K_i-\delta_2V_i\geq 0
\end{align*}
for all \(0\leq i\leq j\), which is equivalent to
\begin{align*}
\frac{\delta_2}{\delta_1}V_i\leq K_i,
\end{align*}
where \(K_i\geq 1\). Making use of the assumptions we obtain the claim.
\end{proof}

\begin{Bem}
\begin{enumerate}
\item A careful analysis of the last proof reveals that 
the condition on the smallness of the scalar potential \eqref{smallness-scalar-potential}
is actually only needed at the singular times \(T_i\), that is
\begin{align*}
\lim_{t\to T_i}\int_M\tilde V(u_t)d\vol_h\leq \frac{\delta_1}{\delta_2}.
\end{align*}
However, it seems rather unlikely that this condition can be satisfied in general.

\item In the case of the standard harmonic map heat flow we have \(\delta_3=1\) and \(V(u_t)=0\)
  such that the bound on the number of singularities reduces to 
   \begin{align*}
    \sum_{i=0}^{j-1}K_i\leq\frac{E(u_0)}{\delta_1}.
   \end{align*}
Moreover, in contrast to the harmonic map heat flow the number of singularities also depends on the metric on \(M\).
\item We want to emphasize that we require the potential \(V(u)\) itself to be sufficiently small and
do not demand any smallness of its gradient. This is what one expects from a mathematics perspective since
the potential itself enters the action functional \eqref{energy-functional}. However, from a physics perspective
the important quantity is the gradient of the potential since it corresponds to the force acting on a system.
\end{enumerate}
\end{Bem}

\begin{Bem}
There is a second way of controlling the number of singularities.
Instead of requiring the potential \(V(u_t)\) to be sufficiently small as in \eqref{smallness-scalar-potential}
we can exploit the fact that \eqref{energy-functional} is not conformally invariant.
More precisely, if we perform a rescaling of the metric \(\tilde h=a h\), where \(a\) is supposed
to be a positive real number, then the first two-terms of \eqref{energy-functional} are not affected,
whereas the scalar potential gets rescaled. More precisely, we find
\begin{align*}
\tilde S_{bos}(\tilde h,u)=\int_M(\frac{1}{2}|du|^2+u^\ast B)d\vol_h+\int_M\frac{1}{a^2}\tilde V(u)d\vol_{\tilde h}.
\end{align*}
In terms of the rescaled metric \(\tilde h=a h\)
the smallness condition \eqref{smallness-scalar-potential}
can be expressed as 
\begin{align*}
\frac{\delta_2}{\delta_1}\int_M\tilde V(u_t)d\vol_{\tilde h}\leq a^2.
\end{align*}
Choosing \(a^2\) large enough we can achieve to have a finite number of singularities
without posing any smallness condition on the scalar potential \(V(\phi)\).
However, the finiteness of the number of singularities now depends on the rescaled metric on the domain \(M\).
\end{Bem}

\subsection{Convergence}
In this subsection we address the issue of convergence of \eqref{flow-rq}.
\begin{Prop}
Let \(u_t\in W\) be a solution of \eqref{flow-rq} and suppose that \(|B|_{L^\infty}<\frac{1}{2}\).
Then there exists a sequence \(t_k\) such that \(u_{t_k}\) converges weakly in $W^{1,2}(M,N)$ 
and strongly in the space $W^{2,2}_{loc}(M\setminus\{x_k,t_k=\infty\})$ 
to a solution of \eqref{euler-lagrange-extrinsic}.
The limiting map \(u_\infty\) is smooth on \(M\setminus\{x_1,\ldots,x_k\}\).
\end{Prop}

\begin{proof}
First, we suppose that $T=\infty$ is non-singular, that is
\[
\limsup_{t\to\infty} (\sup_{x\in M}E(u_t,B_R(x)))<\delta_1
\]
for some \(R>0\).
Since we have a uniform bound on the \(L^2\)-norm of the \(t\) derivative of \(u_t\) 
by Lemma \ref{equality-energy}, we can achieve for $t_k\to\infty$ suitably that
\[
\int_M\big|\frac{\partial u_t}{\partial t}\big|^2d\vol_h\big|_{t=t_k}\to 0.
\]
By \eqref{bound-second-derivatives} we have a bound on the second derivatives
\[
\int_M|\nabla^2u_t|^2(\cdot,t_k)d\vol_h\leq C.
\]
Moreover, by the Rellich-Kondrachov embedding theorem we have
\begin{align*}
u(\cdot,t_k)\to u_\infty& \qquad \textrm{strongly in}~W^{1,p}(M,N)
\end{align*}
for any $p<\infty$.
We get convergence of the evolution equation \eqref{flow-rq} in $L^2$, 
consequently $u_\infty$ is a solution of \eqref{euler-lagrange-extrinsic}
satisfying \(u_\infty\in W^{2,2}(M,N)\).

If $T=\infty$ is singular, that is at the points
$\{x_1,\ldots,x_k\}$
\[
\limsup_{t\to\infty} E(u_t,B_R(x_j))\geq\delta_1,\qquad 1\leq j\leq k
\]
for all $R>0$, then for suitable numbers $t_k\to\infty$ the family $u_{t_k}$
will be bounded in $W^{2,2}_{loc}(M,N)$ on the set \(M\setminus\{x_1,\ldots,x_k\}\).
Consequently, the family \(u_{t_k}\) will accumulate as follows
\begin{align*}
u_\infty\colon M\setminus\{x_1,\ldots,x_k\}\to N.
\end{align*}
We set \(\tilde{M}:=M\setminus\{x_1,\ldots,x_k\}\).
Since we have enough control over the energy of \(u_\infty\) by \eqref{monotonicity-energy-global}, that is
\(E(u_\infty)\leq C\),
we can apply Theorem \ref{theorem-regularity} finishing the proof.
\end{proof}
This completes the proof of Theorem \ref{theorem-flow}.

\begin{Bem}
In the case of a target with negative curvature we have a good understanding of
the properties of the limiting map \(u_\infty\), see \cite[section 4]{MR3624770}.
However, this requires that the second variation of the energy functional is positive,
which makes use of the target having negative curvature.

In the case of the full bosonic string and an arbitrary target most of the methods
employed in the study of harmonic maps can no longer be applied. For this reason
it seems difficult to obtain detailed information on the properties of the limiting map \(u_\infty\)
such that this topic deserves further investigation.
\end{Bem}

\subsection{Blowup analysis}
In order to discuss a blowup analysis of the singular points recall the definition of the 
parabolic cylinder
\begin{align*}
P_r(z_0):=\{z=(x,t)\in M\times (0,\infty) \mid |x-x_0|\leq R, t_0-R^2\leq t\leq t_0\},
\end{align*}
where \(0<R<\min\{\iota_M,\sqrt{t_0}\}\).
Set 
\begin{align*}
v_k(x,t):=u(x_k+r_kx,t_k+r_k^2t),\qquad (x,t)\in P_{r_k^{-1}}.
\end{align*}
For simplicity, assume that \((0,0)\) is a singular point of \(u\in C^\infty(P_1(0,0)\setminus\{0,0\},N)\).
Then there exist \(r_k\to 0\) as \(k\to\infty\) and  \(z_k=(x_k,t_k)\) with \(x_k\to 0, t_k\to 0\) as \(k\to\infty\).
It is easy to check that \(v_k\) satisfies
\begin{align}
\label{vk-equation}
\frac{\partial v_k}{\partial t}=\Delta v_k-\sff(dv_k,dv_k)-Z(dv_k(e_1)\wedge dv_k(e_2))-\frac{1}{r_k^2}\nabla V(v_k).
\end{align}
In the limit \(k\to\infty\), we would have \(P_{r_k^{-1}}\to\R^2\times \R_-\), but it is obvious that
\eqref{vk-equation} blows up as \(k\to\infty\).

This behavior should be expected since the scalar potential \(V(u)\) breaks the conformal invariance 
of the energy functional \eqref{energy-functional}.

However, if \(V(u)=0\) the energy functional \eqref{energy-functional} is invariant
under conformal transformations on the domain and we find that
\begin{align*}
E(u_{t_k},B_{r_k}(x_k))=\sup_{z=(x,t)\in P_1, -1\leq t\leq t_k} E(u_t,B_{r_k}(x))=\frac{\delta_1}{C}
\end{align*}
for \(C>0\) sufficiently large. Assume that \(t_k-4r_k^2\geq -1\).
Moreover, we have
\begin{align*}
\int_{P_{r_k^{-1}}}\big|\frac{\partial v_k}{\partial t}\big|^2d\mu=\int_{t_0-R^2}^{t_0}\int_M\big|\frac{\partial u_t}{\partial t}\big|^2d\vol_h\to 0
\end{align*}
and also
\begin{align*}
E(v_k(t))&\leq E(u_0),\qquad -r_k^{-2}\leq t\leq 0, \\
\sup_{(x,t)\in P_k}E(v_k(t),B_2(x))&\leq C\sup_{(x,t)\in P_1} E(u_t,B_{r_k}(x))\leq\delta_1.
\end{align*}
Consequently, we can take the limit \(k\to \infty\) and \(v_k\) converges to some limiting map \(\omega\).
Then, \(\omega\in C^\infty(\R^2\times (-\infty,0),N)\) solves
\begin{align*}
0=\Delta\omega-\sff(d\omega,d\omega)-Z(d\omega(e_1)\wedge d\omega(e_2))
\end{align*}
since \(\partial_t\omega=0\). Using the conformal invariance we perform a stereographic projection to \(S^2\)
and obtain a solution of
\begin{align*}
\tau(\phi)=Z(\zarg),
\end{align*}
where \(\phi\colon S^2\to N\).
Making use of a Theorem of Grüter \cite{MR744314} we can remove the singular points that 
we get from the stereographic projection. Hence, we get a variant of the usual bubbling
that is well known in the standard harmonic map heat flow.

\begin{Bem}
In the case that \(\dim N=3\) and \(\phi\) is an isometric immersion the
equation
\begin{align*}
\tau(\phi)=Z(d\phi(e_1)\wedge d\phi(e_2))
\end{align*}
is known as \emph{prescribed curvature equation}.
Thus, the \emph{bubbling} described above in the case of \(V(\phi)=0\)
yields maps with prescribed mean curvature from \(S^2\).
However, the condition \(|B|_{L^\infty}<\frac{1}{2}\) that we needed to impose 
does not seem to have a natural geometric interpretation.
\end{Bem}

\subsection{Qualitative properties of the limiting map}
Let us briefly discuss the qualitative behavior of solutions to \eqref{euler-lagrange-intrinsic}.

\begin{Prop}
Let \(\phi\colon M\to N\) be a smooth solution of
\begin{align*}
\tau(\phi)=Z(\zarg)+\nabla V(\phi).
\end{align*}
By \(|\kappa^N|\) we denote an upper bound on the sectional curvature of \(N\). If
\begin{align}
\label{small-energy-assumption}
\frac{\operatorname{Scal}}{2}\geq (|Z|^2_{L^\infty}+|\kappa^N|)|d\phi|^2+|\hess V|_{L^\infty},
\end{align}
then the map \(\phi\) is trivial.
\end{Prop}

\begin{proof}
By a direct calculation we find (see \cite[Lemma 3.1]{MR3573990} for more details)
\begin{align*}
\Delta\frac{1}{2}|d\phi|^2=&|\nabla d\phi|^2+\frac{\operatorname{Scal}}{2}|d\phi|^2-\langle R^N(d\phi(e_\alpha),d\phi(e_\beta))d\phi(e_\alpha),d\phi(e_\beta)\rangle \\
&-\langle Z(\zarg),\tau(\phi)\rangle+\hess V(d\phi,d\phi) \\
\geq& |\nabla d\phi|^2+\frac{\operatorname{Scal}}{2}|d\phi|^2-|\kappa^N||d\phi|^4-|Z|_{L^\infty}|d\phi|^2|\tau(\phi)|-|\hess V|_{L^\infty}|d\phi|^2.
\end{align*}
Using that \(|\tau(\phi)|^2\leq 2|\nabla d\phi|^2\) and applying Young's inequality we deduce
\begin{align*}
\Delta\frac{1}{2}|d\phi|^2&\geq|d\phi|^2(\frac{\operatorname{Scal}}{2}-|Z|^2_{L^\infty}|d\phi|^2-|\kappa^N||d\phi|^2-|\hess V|_{L^\infty})\geq 0,
\end{align*}
where we used the assumptions in the last step.
Consequently, \(|d\phi|^2\) is a subharmonic function and thus has to be constant.
\end{proof}

\begin{Bem}
If we integrate the condition \eqref{small-energy-assumption} over the surface \(M\) we obtain
\begin{align*}
\pi\chi(M)\geq (|Z|^2_{L^\infty}+|\kappa^N|)\int_M|d\phi|^2d\vol_h+|\hess V|_{L^\infty}\vol(M,h).
\end{align*}
Note that this condition can only be satisfied on surfaces of positive genus.
\end{Bem}

\par\medskip
\textbf{Acknowledgements:}
The author gratefully acknowledges the support of the Austrian Science Fund (FWF) 
through the START-Project Y963-N35 of Michael Eichmair and
the project P30749-N35 ``Geometric variational problems from string theory''.
\bibliographystyle{plain}
\bibliography{mybib}
\end{document}